\documentclass[11pt, reqno]{amsart}
\usepackage{graphicx,psfrag}
\usepackage[psamsfonts]{amssymb}
\usepackage{pdfsync}
\usepackage{amscd}
\usepackage{amsmath}
\usepackage{mathrsfs}
\usepackage{stmaryrd}
\usepackage[small,nohug,heads=littlevee]{diagrams}
\usepackage{comment}
\usepackage{enumerate}
\usepackage{xfrac}
\usepackage{cases}
\diagramstyle[labelstyle=\scriptstyle]

\theoremstyle{plain}
    \newtheorem{theorem}{Theorem}[section]
    \newtheorem{lemma}[theorem]{Lemma}
    \newtheorem{corollary}[theorem]{Corollary}
    \newtheorem{proposition}[theorem]{Proposition}
\theoremstyle{definition}
    \newtheorem{definition}[theorem]{Definition}

\theoremstyle{remark}
    \newtheorem{remark}[theorem]{Remark}

\numberwithin{equation}{section}

\newcommand{\T}{\mathbb{T}}

\newcommand{\Z}{\mathbb{Z}}
\newcommand{\C}{\mathbb{C}}
\newcommand{\R}{\mathbb{R}}
\newcommand{\D}{\mathbb{D}}

\newcommand{\To}{\longrightarrow}

\newarrow{Equal} =====

\title{Proper Maps, Bordism, and Geometric Quantization}

\author{Yanli Song}

\date{\today}

\begin{document}

\abovedisplayskip=2pt
\belowdisplayskip=2pt

\maketitle

\begin{abstract}
Let $G$ be a compact connected Lie group acting on a stable complex manifold $M$ thtat has an equivariant vector bundle $E$ on it.  In addition, suppose that $\phi$ is an equivariant map from $M$ to the Lie algebra $\mathfrak{g}$.  We  define an equivalence relation on the triples $(M, E, \phi)$ such that the set of equivalence classes forms an abelian group. We prove that this group is isomorphic to a completion of character ring $R(G)$ and so give a geometric proof of the Quantization Commutes with Reduction conjecture in the non-compact setting. 

\end{abstract}

\tableofcontents

\section{Introduction}
This article focuses on  the Quantization  Commutes with Reduction  conjecture \cite{Guillemin82} in the non-compact setting,  a conjecture proved by Ma and Zhang \cite{Zhang09} (Paradan later gave a different proof \cite{Paradan09}). We provide a new approach to  this conjecture, one that is closely related to noncommutative geometry. In addition, our methods  lead to further extensions. 

In the standard setting of the quantization commutes with reduction conjecture, $M$  is a compact symplectic manifold with symplectic form $\omega$.  Assume that $E$ is a complex line bundle carrying a Hermitian metric and  a Hermitian connection $\nabla^{E}$ such that 

\[ \frac{\sqrt{-1}}{2\pi}(\nabla^{E})^{2} = \omega. \]
Additionally, fix an almost complex structure $J$ such that 

\[ g^{TM}(v, w) = \omega(v, Jw), \ \ v, w \in  TM \]
defines a Riemannian metric on $M$. Let $G$ be a compact connected Lie group, with Lie algebra denoted by $\mathfrak{g}$,  acting on $M$ and $E$ in a Hamiltonian fashion. That is, 

  \[ d\mu(\xi) = \iota_{\xi_{M}} \omega, \]
 where, by the Kostant formula \cite{Kostant70}, the moment map $\mu : M \To \mathfrak{g}^{*}$  is defined by

\begin{equation}
\label{s70}
 \mu(\xi)  = \frac{\sqrt{-1}}{2\pi}(\nabla^{E}_{\xi_{M}} - L_{\xi} ), \ \xi \in \mathfrak{g}.
 \end{equation}
Here, $\xi_{M}$ is the induced infinitesimal vector field and $L_{\xi}$ is the Lie derivative.  In this case, we call $(M, E, \omega, \mu)$ $\emph{pre-quantum  data}$ \cite{Guillemin82}.  Then one can canonically construct a Spin$^{c}$-Dirac operator 
 
 \[  D^{E} : \Omega^{0, *}(M, E) \to \Omega^{0, *}(M, E), \]
 which  gives a finite dimensional virtual vector space
 
 \[ \mathrm{Index}(M, E) = \mathrm{Ker}(D^{E} ) \cap \Omega^{0, \mathrm{even}}(M, E)- \mathrm{Ker}(D^{E} ) \cap \Omega^{0, \mathrm{odd}}(M, E). \]
 
 \begin{definition}
 Given any pre-quantum data $(M, E, \omega, \mu)$, we define its geometric quantization 
 
 \[ Q(M, E) = \mathrm{Index}(M, E),\] 
 which is a virtual representation of $G$.
 \end{definition}

Let $\hat{G}$ be the set of dominant weights in $\mathfrak{g}^{*}$ \cite{Humphreys78}. It is well-known that there is  an one to one correspondence between $\hat{G}$ and all the irreducible G-representations \cite{Kirillov04}. Take any $\gamma \in \hat{G}$. When $\gamma$ is a regular value of the moment map $\mu$ (the singular case was discussed in \cite{Meinrenken99}) , the action of $G$ on $\mu^{-1}(G \cdot \gamma)$ is locally free. Therefore, $M_{\gamma} = \mu^{-1}(G\cdot \gamma) /G$ is an orbifold and $E_{\gamma}= (E|_{\mu^{-1}(G\cdot \gamma)})/G$ is an orbifold line bundle\cite{Weinstein79} (for basic definitions of orbifold, see \cite{Sateka57} \cite{Kawasaki79}). Moreover, $M_{\gamma}$ inherits symplectic structure from $M$. Hence, we can build a Spin$^{c}$-Dirac operator, with $\mathrm{Index}(M_{\gamma}, E_{\gamma})$ given by the orbifold index theorem of Kawasaki \cite{Kawasaki79}. 

\begin{definition}
\label{k3}
We can define the quantization of pre-quantum data $(M, E, \omega, \mu)$ in an alternative way:
 
 \[ Q_{\mathrm{RED}}(M, E, \mu) = \sum_{\gamma \in \hat{G}} \mathrm{Index}(M_{\gamma}, E_{\gamma}) \cdot V_{\gamma} \in R(G), \]
 where $V_{\gamma}$ is the irreducible representation with highest weight $\gamma$.

\end{definition}

 The celebrated quantization commutes with reduction theorem \cite{Guillemin82} says that the two definitions of quantization coincide. 
 \begin{theorem}[Meinrenken-Sjammar, Tian-Zhang, Paradan]
 If $(M, E, \omega, \mu)$ is pre-quantum data, then
 
 \[ Q(M, E) = Q_{\mathrm{RED}}(M, E) \]
 \end{theorem}
 
 Now, we assume that $M$ is noncompact. In order to formulate a suitable $[Q, R]=0$ theorem in noncompact case, the basic problem is ``how to quantize a noncompact manifold".  For the case in which $\mu$ is proper (i.e, the inverse image of a compact subset is compact), Ma and Zhang introduced a formal geometric quantization for $(M, E, \mu)$,  by means of the Atiyah-Patodi-Singer-type index \cite{Atiyah79} for Dirac-type operators on manifolds with a boundary, denoted by $Q_{\mathrm{APS}}(M, E, \mu)$(see \cite{Zhang09} for details).  
 
 On the other hand, since $\mu$ is proper, the reduced orbifold $M_{\gamma} = \mu^{-1}(G\cdot \gamma) /G$ is compact. Thus, we extend  Definition \ref{k3} to the noncompact case\cite{Weitsman01} which we also denote $Q_{\mathrm{RED}}(M, E, \mu)$. It is necessary to point out  that both $Q_{\mathrm{APS}}(M, E, \mu)$ and $Q_{\mathrm{RED}}(M, E, \mu)$ take values in $\hat{R}(G) = \mathrm{Hom}_{\Z}(R(G), \Z)$, the completion of character ring $R(G)$.

 We can now state the result of Ma-Zhang \cite{Zhang09, Paradan09}.

 \begin{theorem}[Ma-Zhang]
 \label{G-S theorem}
 When $(M, E, \mu)$ is pre-quantum data and $\mu$ is proper, we have 
 
 $$ Q_{\mathrm{APS}}(M, E, \mu) = Q_{\mathrm{RED}}(M, E, \mu). $$
 
 \end{theorem}
 
In this paper, we take a topological account of this problem.

 \begin{definition}
A $\emph{stable complex  structure}$ on an orbifold $M$ is an equivalent class of complex structures on $TM \oplus \R^{k}$.  We say two stable complex structures on $TM \oplus \R^{k_{1}}$ and $TM \oplus \R^{k_{2}}$ are equivalent if there exists $r_{1}$ and $r_{2}$ such that 

\[ TM \oplus \R^{k_{1}} \oplus \C^{r_{1}} \ \mathrm{and} \ TM \oplus \R^{k_{2}} \oplus \C^{r_{2}}  \]
are isomorphic complex vector bundles. 
\end{definition}
When a group $G$ acts on M, we define an equivariant stable complex structure by requiring the complex structures to be invariant and the isomorphism to be equivariant. Here, the group acts on $TM$ by the natural lifting of its action on $M$, and it acts trivially on the trivial bundles $\C^{r}$. However, unlike the usual case, we do not require that $G$ acts trivially on $\R^{k}$.

Instead of pre-quantum data, we consider more general  data $(M, E, \phi)$ as follows:

 \begin{itemize}

\item 
   $M$ is  a stable complex $G$-orbifold, possibly noncompact. 
\item
   $E$ is a $G$-equivariant orbifold vector bundle over $M$. 
\item
   $\phi$ is  a $G$-equivariant map from $M$ to $\mathfrak{g} \cong \mathfrak{g}^{*}$ (we can identify $\mathfrak{g}$ and its dual by making a choice of invariant inner product on $\mathfrak{g}$). 
 
 \end{itemize}
 
\begin{definition}
Given any $G$-equivariant map $\phi : M \mapsto \mathfrak{g}$,  define a vector field $V^{\phi}$ by the formula:
  
$$ V^{\phi}(m):= \frac{d}{dt} \Big|_{t=0} \mathrm{exp}(-t\phi(m)) \cdot m, \ \forall m \in M.$$ 
Let $M^{\phi}$ be the $\it{vanishing \ subset}$ in $M$:

$$ M^{\phi} = \{ m \in M \big | V^{\phi}(m) = 0 \}.$$
\end{definition}

We do not require $\phi$ to be moment map. Rather, we relax the moment map condition in the following way. 

\begin{definition}

We say that $\phi$ is \emph{compatible with E},  if there exists a constant $K$ such that  
 \begin{equation}
 \label{s71}
  \| \frac{\sqrt{-1}}{2\pi}L_{\xi}  +  \langle \phi(m), \xi \rangle \cdot I_{m}   \|  \leq K \| \xi\|, \ \text{for \  all}  \ m \in M^{\phi} ,
 \end{equation}
  where  $I_{m}$ is the identity map from $E|_{m}$ to itself, $L_{\xi}$ is the Lie derivative on $E|_{m}$,  and $\xi$ lies in the isotropy Lie algebra $\mathfrak{g}_{m}$. 
 
\end{definition}

\begin{remark}
When $E$ is an actual orbifold line bundle with a moment map $\mu$,  (\ref{s71}) is equivalent to 

\begin{equation}
\label{s27}
| \langle \mu(m), \xi \rangle - \langle \phi(m) , \xi \rangle | \leq K \cdot \| \xi \|.
\end{equation}
Notice that (\ref{s27}) does not depend on the choice of moment map $\mu$.

\end{remark}

 \begin{definition}
We say that a triple $(M, E, \phi)$ in which $M$ may have boundary is a $\emph{K-chain}$ if  $\phi$ is proper over $M^{\phi}$ and  compatible with $E$.  When $M$ has no boundary, we say that $(M, E, \phi)$ is a $\emph{K-cycle}$. 

\end{definition}
In particular, all  pre-quantum data $(M, E, \mu)$ and compact triples are $K$-cycles.


 Next, we will define the equivalence relation between $K$-cycles. To begin with,   a bordism between two $n$-dimensional orbifolds $M$ and $N$ is a $n+1$ dimensional orbifold $W$ with boundary, such that 
 \[  \partial W = M \sqcup (-N). \]
In ordinary bordism theory, we only consider compact orbifold (otherwise, any orbifold $M$ is bordant to the empty set via $W = [0, 1) \times M$). In order to obtain a nontrivial theory, we use the following definition. 
 \begin{definition} 
 
 \label{s12}
Suppose that $(W, L, \psi)$ is a $K$-chain, whose boundary is divided into two parts $ M \sqcup N$. Hence, we obtain two $K$-cycles:

\[ (M, L|_{M}, \psi|_{M}) \ \mathrm{and} \ (N, L|_{N}, \psi|_{N}). \]
We say that the first $K$-cycle $\it bordant$ to (the opposite of) the second. 
 
 \end{definition}
 
 \begin{definition}
 \label{d1}
 
Our equivalence relation between $K$-cycles is generated by the following three elementary steps:

\begin{description}
\item[Disjoint\ Union]
 $ (M, E, \phi) \bigsqcup (M, F, \phi) \sim (M, E \oplus F, \phi). $
    
\item[Bundle\ Modification   \cite{Baum07}]
 Suppose that $P$ is a principal bundle over $M$ whose structure group is the compact Lie group $H$. Let $N$ be a compact, even dimensional, stable complex $H$-orbifold.  The orbifold index of the associated Spin$^{c}$-Dirac operator gives an element $[D_{N}] \in R(H)$. And, if $[D_{H}] = [1]$, then

\[ (M, E, \phi) \sim (\hat{M}, \hat{E}, \hat{\phi}), \]
 where $\hat{M} = P \times_{H} N$, $\hat{E}$ is the pull back of $E$ and $\hat{\phi}$ is the composition of $\phi$ with the projection to $M$. 
   
\item[Bordism]
  Definition \ref{s12}.
     
\end{description}

\end{definition}

\begin{definition}

\label{def-K}

We denote by $\hat{K}(G)$ the set of  equivalence classes of $K$-cycles. 

\end{definition}

The set $\hat{K}(G)$ is an abelian group, whose addition operation is given by disjoint union; the additive inverse of a $K$-cycle is obtained by reversing the stable complex structure. 

For any $K$-cycle $(M, E, \phi)$, we use $n \cdot (M, E, \phi)$ to denote the disjoint union of $n$ copies of $(M, E, \phi)$.  Suppose $\Gamma = \sum_{\gamma \in \hat{G}} n_{\gamma} V_{\gamma}$ is an arbitrary element in $\hat{R}(G)$.  Define

\[ \mathcal{O}_{\Gamma} = \bigsqcup_{\gamma \in \hat{G}}n_{\gamma} \cdot (\mathcal{O}_{\gamma}, E_{\gamma}, \iota_{\gamma}), \]
where $\mathcal{O}_{\gamma}$ is the coadjoint orbit through $\gamma \in \mathfrak{t}_{+}$; $E_{\gamma}$ is the natural line bundle defined using weight $\gamma$ \cite{Bott65} \cite{Kirillov04}; and $\iota_{\gamma}$ is the inclusion.  It is clear that $ \mathcal{O}_{\Gamma}$ is a $K$-cycle. 

The following theorems constitute the main results of this paper. 

\begin{theorem}The map $P : \hat{R}(G) \to\hat{K}(G)$
\label{o9}

\[ P  : \Gamma \longmapsto \mathcal{O}_{\Gamma}\] 
gives an isomorphism of abelian groups and $R(G)$-modules. 
\end{theorem}

 The idea of the definition of $\hat{K}(G)$ comes from  the geometric $K$-homology defined by Baum and Douglas \cite{Baum81}. We can generalize $\hat{K}(G) = \hat{K}^{G}(\mathrm{pt})$ so as to obtain not just a group but a functor:
 
 \[ X \mapsto \hat{K}^{G}(X), \]
 where $X$ is a paracompact Hausdorff $G$-space. To be  more precise, let $(M, E, \phi, f)$ be a 4-tuple, where $M, E, \phi$ are the same as in the definition of $K$-cycles  and $f$ is an equivariant map from $M$ to $X$. The equivalence relation extends naturally to this general case. In a work that is nearly finished, we show that the set of equivalence classes $\hat{K}(G)(X)$ is isomorphic to the Kasparov group $KK(C^{*}(G,X), \C)$ \cite{Kasparov88}. 

\begin{theorem}
\label{main theorem}
The inverse map $Q_{\mathrm{TOP}} = P^{-1} : \hat{K}(G) \To \hat{R}(G)$ has the following properties:

\begin{enumerate}
 
 \item
 When $M$ is compact,  $ Q_{\mathrm{TOP}}(M, E, \phi) = Q(M, E) \in R(G).$
 
 \item
When $(M,E, \phi_{1})$ and $(N, F, \phi_{2})$ are two K-cycles and $N$ is compact, we have 

\[ Q_{\mathrm{TOP}}(M, E, \phi_{1}) \times Q_{\mathrm{TOP}}(N, F, \phi_{2}) = Q_{\mathrm{TOP}}(M \times N, E \boxtimes F, \hat{\phi}_{1} + \hat{\phi}_{2}), \]
where $\hat{\phi}_{1}$  and  $\hat{\phi}_{2}$ are the pullbacks of $\phi_{1}$ and $\phi_{2}$. 

\end{enumerate}

\end{theorem}

According to (1),  $Q_{\mathrm{TOP}}$ can be considered as a generalization of the usual quantization for compact manifold. According to (2), $Q_{\mathrm{TOP}}$ satisfies the ``multiplicative property", which is one of the main difficulties in \cite{Zhang09}. 

\begin{theorem}
\label{o7}
 When $(M, E, \phi)$ is pre-quantum  data and $\phi$ is proper, we have
 
 \[ Q_{\mathrm{TOP}}(M, E, \phi) = Q_{\mathrm{RED}}(M, E, \phi). \]
 
 \end{theorem}

$\mathbf{Acknowledgments}$: I am grateful to my advisor, Professor N. Higson, for his kind guidance and advice. Also, I would like to thank Professor W. Zhang and Professor P. Baum for the many helpful conversations we have shared.



\section{Basic Properties of $K$-Cycles}
In this section, we will discuss some of the basic properties of $K$-cycles. 

\begin{definition}
For convenience, we define some special $K$-chain ($K$-cycle):
\begin{itemize}

\item
  We say that a $K$-chain or a $K$-cycle $(M, E, \phi)$  $\it has \ compact \ vanishing \ set$ if $M^{\phi} \subseteq M$ is compact. 
  
\item
  We say that a $K$-cycle is $\it closed$ if $M$ is  compact and has no boundary. 
   
\item
  We say that a $K$-cycle is $\it discrete$  if it has the following form:

\[ \bigsqcup_{k=1}^{\infty}(N_{k}, E_{k}, \rho_{k}) \in\hat{K}(G), \]
where $N_{k}$ are closed orbifolds. 

\end{itemize}
\end{definition}

These lemmas follow immediately from the definition of $K$-cycles.

\begin{lemma}
 When $M$ is compact, we have  
 
\[ (M, E, \phi) \sim (M, E, 0). \]
\end{lemma}

\begin{lemma}
Every  $K$-cycle $(M, E, \phi)$ is equivalent to a finite sum of $K$-cycles 

\[ (M, E, \phi) \sim \sum_{i=1}^{n} (M_{i}, E_{i}, \phi_{i}) \]
where $\{ E_{i}\}$ are line bundles. 
\begin{proof}
By bundle modification, we have 

\begin{equation}
(M, E, \phi) \sim (\hat{M}, \hat{E}, \hat{\phi}),
\end{equation}
where $\hat{M}$ is a bundle of flag manifold over $M$ induced by $E$; $\hat{E}$ is the pullback of $E$; and $\hat{\phi}$ is the composition of $\phi$ with the projection from $\hat{M}$ to $M$. Hence, $\hat{E}$ splits into sum of line bundles. 

\end{proof}
\end{lemma}

 Let us recall the definition of the geometric $K$-homology by Baum and Douglas \cite{Baum81} \cite{Baum10}. Let $X$ be a paracompact Hausdorff $G$-space. A cycle is a triple  $(M, E, f)$ consisting of:

\begin{itemize}
\item
  $M$ = a weakly complex compact manifold with $G$-action. 
  \item
  $E$ =  a $G$-equivariant vector bundle over $M$.
  \item
   $f$ = a G-equivariant continuous map from $M$ to $X$. 
\end{itemize}
 The equivalence relations are generated by direct sum, bordism (in the compact sense) and bundle modification. The equivalence classes of cycles form an abelian group $K^{G}(X)$. In particular, when $X$ is a point, the map $f$ is trivial. Therefore, we have a natural map:
 
 \[ B : K^{G}(\mathrm{pt}) \To\hat{K}(G) : B(M, E) = (M, E, 0).\]

 \begin{theorem}[Localization Theorem]

\label{localization theorem}
Suppose $\{ U_{\alpha} \}$ is a family of disjoint G-invariant open subsets such that 

\[ M^{\phi} \subseteq \bigsqcup_{\alpha} U_{\alpha}. \]
We have 

$$  (M, E, \phi) \sim \bigsqcup_{\alpha}(U_{\alpha}, E|_{U_{\alpha}}, \phi|_{U_{\alpha}}). $$

\begin{proof}
Let $W = M \times [0, 1]$ and $\hat{E}$ be the pullback of $E$. Define a map 

$$\hat{\phi} : W \To \mathfrak{g} :(m, t) \rightarrow \phi(m).$$ 
In addition, let $F = M \setminus (\bigcup(U_{\alpha})$ and $\hat{W} = W \setminus (F \times \{ 1 \})$. We also denote by $\hat{E}$ and $\hat{\phi}$ their restrictions to $\hat{W}$.  It is easy to verify that 

\[ (\hat{W}, \hat{E}, \hat{\phi}) \]
is a $K$-chain which gives desired bordism.

\end{proof}

\end{theorem}

\begin{remark}

If $G$ is the circle group $S^{1}$,  Theorem \ref{localization theorem} is similar to the Linearization Theorem in \cite{Ginzburg96} in which $\phi$ is required to be an abstract moment map. 
\end{remark}

\begin{corollary}[Gluing Property]
\label{gluing}

Let  $(M, E, \phi)$ be any $K$-cycle. Suppose that $\Sigma \subset M$ is a smooth G-invariant hypersurface in M and $\Sigma$ cuts M into  two oriented piece: $M\setminus \Sigma= M_{+}\sqcup M_{-}$. We obtain two $K$-cycles: 

\[ (M_{+}, E|_{M_{+}}, \phi_{+}) \ \mathrm{and} \ (M_{-}, E|_{M_{-}}, \phi_{-}). \] 
When the vector field $V^{\phi}$ does not vanish  over $\Sigma$, we have 

$$(M, E, \phi) \sim (M_{+}, E|_{M_{+}}, \phi_{+}) + (M_{-}, E|_{M_{-}}, \phi_{-}). $$
\end{corollary}



\section{Vanishing set of $K$-Cycles}
The main goal of this section and the next one is to study $K$-cycles with compact vanishing. These two sections are the main part of this paper. 

For any arbitrary $K$-cycle $(M, E, \phi)$,  the vanishing set $M^{\phi}$ will play an important role in defining its quantization. In general, $M^{\phi}$ is a subset of $M$, which may be very complicated. However,  the following theorem shows that $M^{\phi}$ can be separated into compact parts. 
 
\begin{theorem}

\label{separate theorem}
Let $(M, E, \phi)$ be a $K$-cycle in which $E$ is a line bundle. There exists a covering of $M^{\phi}$ by $G$-invariant disjoint open subsets of $M$, $\{ U_{\alpha} \}$ such that each $F_{\alpha} = U_{\alpha} \cap M^{\phi}$ is compact. 

\end{theorem}

\begin{remark}

When $(M, E, \phi)$ is pre-quantum data and $\phi$ is proper,  Theorem \ref{separate theorem} is trivial.  In fact, let $\mathcal{H} = \| \phi \|^{2} : M \To \R$. We can find a series of regular values of $\mathcal{H}$ 

\[ c_{1}, c_{2}, \dots, c_{n}, \dots \] 
such that  $\lim_{n \to \infty}c_{n} =  \infty$.  Over $\mathcal{H}^{-1}(c_{i}) \subseteq M$, the vector field $V^{\phi}$ does not vanish.  Hence,

\[ \{ \ U_{i} = \mathcal{H}^{-1}((c_{i} , c_{i+1})) \}_{i=1}^{\infty} \ \]  
give desired covering. 
\end{remark}

 In general, the proof of Theorem \ref{separate theorem} is based on Lemma \ref{y4} and Lemma \ref{y8}. To begin with,  let  $\T$ be  a maximal torus in $G$, $\mathfrak{t}$ be the Lie algebra of $\T$, and $\mathfrak{t}_{+}$ be a chosen positive Weyl chamber. We observe that $m \in M^{\phi} $ if and only if $m \in M^{\gamma} \cap \phi^{-1}(\gamma)$, where $\gamma = \phi(m) \in \mathfrak{g}$. Thus, 

\begin{equation} 
\label{y3}
M^{\phi} = \bigcup_{\gamma \in \mathfrak{g}}(M^{\gamma} \cap \phi^{-1}(\gamma)) = \bigcup_{\gamma \in \mathfrak{t}_{+}}G.(M^{\gamma} \cap \phi^{-1}(\gamma)). 
\end{equation}


Consider  the set  of stabilizers $\{ H_{i} \}_{i=1}^{\infty}$ for the  action of maximal torus $\T$ on  $M$. Put an invariant connection on $E$ and denote by $\mu$ the associated moment map. Since $H_{i}$ are subgroups of $\T$, 
we can identify their Lie algebras $\mathfrak{h}_{i}$ as subspaces in $\mathfrak{t}$.  Let 

\[ \mu_{\mathfrak{t}} : M \to \mathfrak{t}\]
be the moment map associated to the action of $\T$ and $\mu_{\mathfrak{h}_{i}}$  be  the composition of $\mu_{\mathfrak{t}}$ with the projection  from $\mathfrak{t}$ to $\mathfrak{h}_{i}$. From (\ref{y3}), we have

\begin{equation}
\label{y9}
 M^{\phi} = \bigcup_{\mathfrak{h}_{i} \in \Gamma}G\cdot (M^{H_{i}} \cap \phi^{-1}(\mathfrak{h}_{i})),
\end{equation}
where $\Gamma$ is the set of Lie algebras $\mathfrak{h}_{i}$ such that $\mathfrak{h}_{i} \cap \mathfrak{t}_{+} \neq \{ 0\}$. For the decomposition in  (\ref{y9}), we have the following lemma.

\begin{lemma}
\label{y4}
For every $\mathfrak{h}_{i} \in \Gamma$, the set 

\[G \cdot (M^{H_{i}} \cap \phi^{-1}(\mathfrak{h}_{i}))\]
can be separated by   $G$-invariant disjoint open subsets $\{ V^{i}_{k} \}_{k=1}^{\infty}$ in $(G\cdot M^{H_{i}})$ such that $\phi$ is bounded on all $V^{i}_{k}$. 

\begin{proof}
Notice that  the map $\mu_{\mathfrak{h}_{i}}$ is locally constant on $M^{H_{i}}$. The set
  
 \[  \mu_{\mathfrak{h}_{i}}(M^{H_{i}})  \cap \mathfrak{t}_{+}\] 
 is a set of countable many points in $\mathfrak{h_{i}} \cap \mathfrak{t}_{+}$. Hence, we can find a covering of $G \cdot (M^{H_{i}} \cap \phi^{-1}(\mathfrak{h}_{i}))$, consisting of disjoint open subsets $\{ V^{i}_{k} \}_{k=1}^{\infty}$ in $(G\cdot M^{H_{i}})$ such that 

\begin{equation}
\label{y2} 
  \mu_{\mathfrak{h}_{i}}(V^{i}_{k}) \cap \mathfrak{t}_{+} = \gamma_{k},
\end{equation}
 where $\gamma_{k} \in \mathfrak{h}_{i} \cap \mathfrak{t}_{+}$. In addition, we can choose $V^{i}_{k}$ small enough such  that 
 
 \begin{equation} 
\label{s10}
\phi(V^{i}_{k}) \cap \mathfrak{t}_{+} \subseteq \{ x \in \mathfrak{t}_{+} \big| \ \text{distance}(x, \mathfrak{h}_{i}) < \epsilon \},
\end{equation}
where $\epsilon$ is arbitrary small positive number. From (\ref{s27}), (\ref{y2}), and (\ref{s10}),  we know that 
 
\[   \phi(V^{i}_{k}) \cap \mathfrak{t}_{+} \subseteq B(\gamma_{k}, K + \epsilon)= \{x \in \mathfrak{t}_{+}\big| \|x - \gamma_{k}\| < K + \epsilon \}. \]
 It follows that $\phi$ is bounded on all $V_{k}^{i}$. 
\end{proof}

\end{lemma}

Obviously,  the set $\bigcup_{i, k}V^{i}_{k}$  covers  the vanishing set $M^{\phi}$. We define $\{U_{\alpha} \}$ to be disjoint $G$-invariant neighborhoods of connected components  of $\bigcup_{i, k}V^{i}_{k}$ in $M$. Thus, the open sets $\{ U_{\alpha} \}$ give a covering of $M^{\phi}$.   It remains to show that $F_{\alpha} = M^{\phi} \cap U_{\alpha}$ is compact. 

\begin{lemma}
\label{y8}
The set $F_{\alpha} = M^{\phi} \cap U_{\alpha}$ is compact for all $\alpha$.  

\begin{proof}
Because $\phi$ is proper over $M^{\phi}$, it is equivalent to prove that $\phi$ is bounded on $U_{\alpha}$. 

Suppose $\phi$ is unbounded on $U_{\alpha}$. According to Lemma \ref{y4}, there must exist  an infinite chain in $\{ V_{k}^{i} \}$,

\begin{equation} 
\label{s5}
V_{1}, V_{2}, \dots, V_{m} \dots 
\end{equation}
such that 

\[ V_{i} \cap V_{i+1} \neq \emptyset \ \ \text{and} \  \ V_{i} \subseteq U_{\alpha}. \] 
Notice that every $V_{i}$ is contained in some $(G\cdot M^{H_{k}})$. From (\ref{y2}), let us denote 

\begin{equation}
\label{y7}
\gamma_{i} = \mu_{\mathfrak{h}_{k}}(V_{i}) \cap \mathfrak{t}_{+}. 
\end{equation}
It is clear that 

\begin{equation}
\label{s11}
\lim_{i \to \infty} \| \gamma_{i} \| = \infty.
\end{equation}
Without loss of generality, we can assume that for any constant $T$, there exists $N_{T} \geq T$ and a chain

\[ \gamma_{N_{T}}, \dots, \gamma_{N_{T}+\text{dim} (\T)},  \]
such that they are two by two distinct.

 For $N_{T} \leq k \leq N_{T}+\text{dim} (\T)-1$, pick an arbitrary point $m_{k} \in V_{k} \cap V_{k+1} \cap \phi^{-1}(\mathfrak{t}_{+})$ and denote by $H_{m_{k}}$ the isotropy group of  $\T$-action. Let $\omega_{k}$ be the $H_{m_{k}}$-weight of the $E|_{m_{k}}$. Here, we can identify $\omega_{k}$ with a point in $\mathfrak{t}$. By (\ref{y7}), we have that
\begin{equation}
\label{s13}
 \omega_{k} -\gamma_{k} \perp \gamma_{k} \ \text{and} \   \omega_{k} -\gamma_{k+1} \perp \gamma_{k+1}.
\end{equation}
In addition, from (\ref{s10}), we have 

\begin{equation}
\label{y5}
\mathrm{distance}(\phi(m_{k}), \mathfrak{h}_{k}) < \epsilon.
\end{equation}
It follows from (\ref{s27}) that

\begin{equation}
\label{y6}
 \|\omega_{k} - \phi(m_{k})\| \leq K.
 \end{equation}
 From (\ref{s13})-(\ref{y6}), we get 
 
 \[ \| \gamma_{k} - \omega_{k} \| \leq K+ \epsilon. \]
 Similarly, we also have 
 
  \[ \| \gamma_{k+1} - \omega_{k} \| \leq K+ \epsilon. \]
Therefore, these points $\{ \omega_{k}, \gamma_{k} \}_{k=N_{T}}^{N_{T} + \text{dim} (\T)}$ are within finite distance from each other. Moreover, they must lie in  the integer lattice of $\mathfrak{t}$ since they are all weights. By the orthogonal condition (\ref{s13}),  there are only finitely many  possibilities of $\{ \omega_{k}, \gamma_{k} \}_{k=N_{T}}^{N_{T} + \text{dim} (\T)}$. This leads to a contradiction to $(\ref{s11})$. 

\end{proof}
\end{lemma}

 This complete the proof of Theorem \ref{separate theorem}. Additionaly, from the proof, one can see  that the assumption that $E$ is a line bundle is not necessary.

\begin{remark}
\label{s30}
In Theorem \ref{separate theorem}, the cover $\{ U_{\alpha} \}$ has the property  that $\phi$ is uniform bounded. That is, there exists a constant $R$ such that for all $\alpha$, 

\[ \| \phi(x) - \phi(y) \| \leq R, \ \mathrm{for} \ x, y \in U_{\alpha} \cap \phi^{-1}(\mathfrak{t}_{+}). \]
\end{remark}

\begin{remark}
\label{regular}
Suppose that $(M, E, \phi)$ is a $K$-cycle with compact vanishing set. Let $\partial M = \overline{M} \setminus M$.  Given any small neighborhood $U$ of $\partial M$ in $\overline{M}$, we can identify 

\begin{equation}
\label{s65}
U \cong \partial M \times [0,1). 
\end{equation}
In fact, by rescaling, we can assume that $\phi(m)$ tends to infinity as $m$ tends to $\partial M$. Take $\mathcal{H} = \| \phi \|^{2} : M \to \R$ and pick a  regular value $c$. Let $M_{c}$ be a  subset of $M$ defined by  

\[  M_{c} = \{ m \in M \big| \mathcal{H}(m) < c \}. \]
When $c$ is large enough, we have $M^{\phi} \subseteq M_{c}$. By Theorem \ref{localization theorem}, 

\[ (M, E, \phi) \sim (M_{c}, E|_{M_{c}}, \phi|_{M_{c}} ),\] 
where the second $K$-cycle satisfies (\ref{s65}). Unless stated otherwise, from here when we refer to the $K$-cycle with compact vanishing set, we always assume that it automatically satisfies $(\ref{s65})$.   In this case, we can extend the vector bundle E and $\phi$ to $\partial M$, denoted by $\partial E$ and $\partial \phi$ respectively. In addition, we denote by $\partial(M, E, \phi)$ the $K$-cycle $(\partial M, \partial E, \partial \phi)$. 

\end{remark}

\begin{remark}
\label{o1}
Suppose $(M_{1}, E_{1}, \phi_{1})$ and $(M_{2}, E_{2}, \phi_{2})$ are two $K$-cycles with compact vanishing set. Assume that there is a  diffeomorphism  $f : \partial M_{1} \cong \partial M_{2}$. By Remark \ref{regular}, the map 
$f$ also induces a diffeomorphism :

\[ \hat{f} : U_{1} \cong U_{2},\]
where $U_{i}$ are neighborhoods of $\partial M_{i}$ as in (\ref{s65}). When the map $\hat{f}$ lifts to an isomorphism between vector bundles $E_{1}|_{U_{1}}$ and $E_{2}|_{U_{2}}$, we can obtain a compact $K$-cycle by gluing the two $K$-cycles using the map $\hat{f}$. In the gluing process, we do not require that $\partial \phi_{1} = \partial \phi_{2}$ because we can alway vary the map $\phi$ without changing $K$-cycle class. 

\end{remark}

  

\section{$K$-Cycles with Compact Vanishing Set}
From the previous two sections,  we know that it is enough to study the $K$-cycles with compact vanishing set. Suppose that $(M, E, \phi)$ is  a $K$-cycles with compact vanishing set. The general strategy is to  build a ``cap'', another $K$-cycle, so that we can compactify $(M, E, \phi)$ by gluing on the cap. The geometric  construction of the cap is the main part of this section. To be more prices, we are going to the prove the following theorems.

\begin{theorem}
\label{a21}
Let $(M, E, \phi)$ be a $K$-cycle with compact vanishing set. There is a $K$-cycle with compact vanishing set $(W, L, \psi)$ such that 

\[ \partial(W, L, \psi) \cong \partial(M, E, \phi). \]
In addition, the $K$-cycle $(W, L, \psi)$ is bordant to a discrete $K$-cycle. 
\end{theorem}

\begin{corollary}
\label{compact theorem}
Every $K$-cycle with compact vanishing set is bordant to a discrete $K$-cycle. 
\end{corollary}

\begin{corollary}
Suppose that $\Sigma$ is a closed $G$-manifold and $\phi : \Sigma \to \mathfrak{g}$  is an equivariant map. If the vector field $V^{\phi}$ induced by $\phi$ is nowhere vanishing over $\Sigma$, then $\Sigma$ is a boundary. 
\end{corollary}


\subsection{Circle Case}
Let $(M, E, \phi)$ be a $K$-cycle with compact vanishing set. According to our discussion in Remark \ref{regular}, we can assume that $M$ is  interior of some compact orbifold with boundary $\Sigma = \partial M$ on which $G = S^{1}$ acts locally freely. 

Suppose that $D^{2}$ is the open disk with  standard $S^{1}$-action.  We can define an orbifold by 

\[  W  = \Sigma\times_{S^{1}} D^{2}, \]
and an orbifold vector bundle on $W$ by 

\[ L = \pi^{*}(\partial E)/S^{1}, \] 
where $\pi$ is the projection of $\Sigma \times D^{2}$ to $\Sigma$, and $S^{1}$-action is the diagonal action.  In addition, the map $(\partial \phi) \circ \pi$ over $\Sigma \times D^{2}$ descends to a map over $W$ denoted by $\psi$. We can verify that $(W, L, \psi)$ indeed constitute a $K$-cycle with compact vanishing set such that 

\[  \partial(W, L, \psi) \cong \partial (M, E, \phi).  \]
As in Remark \ref{o1}, we can obtain a compact $K$-cycle by  gluing  $(W, L, \psi)$ and $(M, E, \phi)$ together.  

Next, we will show that $(W, L, \psi)$ is bordant to a discrete $K$-cycle, beginning with the following lemma.

\begin{lemma}
\label{lemma-2}
Consider the $K$-cycle

\[(D^{2}, \C, f(z)),\]
where  $\C$ is the trivial line bundle over $D^{2}$ on which $S^{1}$ acts trivially, and $f(z)$ is a positive function on $D^{2}$. Then, we have that

\begin{equation}
\label{a1}
(D^{2}, \C, f) \sim \sum_{n=0}^{\infty} (S^{2}, F_{n}, f_{n}), 
\end{equation}
where $S^{2}$ is the sphere with standard $S^{1}$-action,  $F_{n}$ are trivial line bundles over $S^{2}$ on which $S^{1}$ acts with $S^{1}$-weight $n$, and $f_{n}$ are equivariant  functions on $S^{2}$ such that $f_{n}$ equals to $n$ when we restrict to the south pole and north pole.
\begin{proof}
Without losing generality,  we  assume that $f_{n}$ is positive along the equator.  Thus, we can break the compact $K$-cycle $(S^{2}, F_{n}, f_{n})$ into two $K$-cycles with compact vanishing set: 

\[ (S^{2}, F_{n}, f_{n}) \sim (S_{+}, F_{n}|_{S_{+}}, f_{n}|_{S_{+}}) + (S_{-}, F_{n}|_{S_{-}}, f_{n}|_{S_{-}}),\]
where $S_{\pm}$ are the hemispheres. In particular, we have that 

\[ (S_{+}, F_{0}|_{S_{+}}, f_{0}) \cong (D^{2}, \C, f).\] 
Next, let $(S^{2}, F_{n}^{n+1})$ be the compact $K$-cycle obtained by gluing 

\[ (S_{-}, F_{n}|_{S_{-}}, f_{n}) \ \mathrm{and} \ (S_{+}, F_{n+1}|_{S_{+}}, f_{n+1}). \]
Here, $F_{n}^{n+1}$ is an equivariant line bundle on $S^{2}$ with fiber weights equal to $n$ at south pole and $n+1$ at north pole. By Atiyah-Bott fixed point theorem, we have that 

\begin{equation}
\label{a10}
\mathrm{Index}(S^{2}, F_{n}^{n+1}) = 0 \in R(S^{1}).
\end{equation}
Hence, according to  \cite{Baum07}, we can conclude that the $K$-cycle $(S^{2}, F_{n}^{n+1})$ is equivalent to an empty $K$-cycle. This completes the proof. 

\end{proof}

\end{lemma} 

\begin{remark}
If $f(z)$ is a negative function on $D^{2}$, then we have a similar result:
\begin{equation}
\label{a2}
(D^{2}, \C, f) \sim -\sum_{n=1}^{\infty} (S^{2}, F_{-n}, -f_{n}).
\end{equation}

\end{remark}

\begin{proposition} 
\label{p1}
The $K$-cycle $(W, L, \psi)$ is bordant to a discrete  K-cycle. 
\begin{proof}
Let $P = \Sigma \times_{S^{1}} S^{2}$. An orbifold vector bundle $L_{n}$ over $P$ is defined by 

\[ L_{n} = [ E \boxtimes F_{n}]/S^{1},\]
where $S^{1}$-action is the diagonal action.  Additionally, the function 

\[ \Psi : \Sigma \times S^{2} \to \R : \Psi(m, x) = \phi(m) + f_{n}(x) \]
descends to a function on $P$, denoted by $\psi_{n}$. For each $n$, we obtain a compact $K$-cycle $(P, L_{n}, \psi_{n})$.  Using Lemma \ref{lemma-2}, one can show that if $\psi$ is positive, then

\[  (W, L, \psi) \sim \sum_{n=0}^{\infty} (P, L_{n}, \psi_{n}). \]
And if $\psi$ is negative, then we have that 

\[  (W, L, \psi) \sim -\sum_{n=-1}^{-\infty} (P, L_{n}, \psi_{n}). \]

\end{proof}

\end{proposition}


\subsection{Torus Case}
Now, let us assume that  $G = \T$ is a torus  and $\Sigma$ is a compact orbifold together with $\T$-action. 

\begin{definition}
\label{o2}
Let $\{ U_{i} \}_{i=1}^{n}$ be an open cover of $\Sigma$. We say $\{ U_{i}, S_{i} \}_{i=1}^{n}$ is a $\emph{good cover}$ if 

\begin{itemize}
\item
    Every $U_{i}$ is $\T$-invariant. 
\item    
    The circle action $S_{i}$ is a factor in some presentation $\T = S^{1} \times \dots \times S^{1}$ and $S_{i}$ acts locally freely on $U_{i}$. 
\item
 Let us define 
\[ \mathcal{A} = \{ I \subseteq \{ 1, \dots, n\} \big| \ U_{I} = \bigcap_{i \in I} U_{i} \neq \emptyset \}. \]
For all $I \in \mathcal{A}$, $\{ S_{i} \}_{i \in I}$ generate an $|I|$-dimensional subgroup $S_{I}$ and it acts locally freely on $U_{I}$. 
   
   \end{itemize}

\end{definition}

Suppose that $\Sigma$ is a compact orbifold with a good cover $\{U_{i}, S_{i}\}_{i=1}^{n}$. We want to construct a cap for $\Sigma$. If we naively build the caps locally, that is, defining $W_{i} = U_{i} \times _{S_{i}} D^{2}$ as in the circle case, then there is a problem that $\{ W_{i} \}$ may not be glued together. In order to overcome this difficulty, we define the local caps in a more subtle way using the compatibility conditions in Definition \ref{o2}. This idea comes from cutting surgery by Lerman \cite{Lerman95}. 

\begin{lemma}
\label{a14}
For any $I \in \mathcal{A}$, we can define an orbifold $W_{I}$ with boundary isomorphic to $U_{I}$. Moreover, for any $I, J \in \mathcal{A}$ and $I \subsetneq J$, there exists a diffeomorphism   $\Phi_{I}^{J}$ from an open subset $U_{I, J}$ of $W_{I}$ to an open subset  $U_{J, I}$ of  $W_{J}$: 

\[ \Phi_{I}^{J} : U_{I, J} \cong U_{J, I} \]
with the property that for any $K \in \mathcal{A}$ and $J \subsetneq K$, 

\begin{equation}
\label{a13}
U_{I, K} = U_{I, J} \cap (\Phi_{I}^{J})^{-1}(U_{J, K}) \ \mathrm{and} \ \Phi_{I}^{J} \circ \Phi_{J}^{K} = \Phi_{I}^{K}.
\end{equation}

\begin{proof}
To begin with, we choose a partition of unity $\{ \varphi_{i} \}_{i=1}^{n}$  subordinate to the open cover $\{ U_{i} \}_{i=1}^{n}$. Let us fix $n$ numbers  $\{ \alpha_{1}, \dots, \alpha_{n} \}$  such that 
\begin{itemize}
\item
   $\alpha_{i} > 1$ for all i, and
\item
   For any $I \in \mathcal{A}$,  $\prod_{i \in I}\alpha_{i}$ is a regular value for map
\begin{equation}
\label{a17}
  \rho_{I} : U_{I} \times [1, \infty)  \to \R^{|I|}  : (\rho_{I,i})_{i \in I} =  ( t \cdot \varphi_{i}(m))_{i\in I}.
\end{equation}
\end{itemize}

For any $I \in \mathcal{A}$, let us define 

\begin{equation}
\label{a12}
  \check{I} = \bigcup_{K \in \mathcal{A}, I \subseteq K} (K \setminus I).  
\end{equation}
Let $V_{I}$ be an open subset in $U_{I} \times [1, \infty)$ defined by

\[ V_{I} = \{ (m, t) \in U_{I} \times [1, \infty) \big| t \cdot \varphi_{i}(m) < \alpha_{i}, \mathrm{for \ all} \ i \in \check{I} \}. \]
Apparently, $V_{I}$ has a boundary isomorphic to $U_{I}$ and $S_{I}$ acts locally freely on $V_{I}$. 

Suppose that $\C^{|I|}$ is the product of $|I|$ copies of $\C$, where every copy  has a individual standard circle action. This gives an $|I|$-dimensional torus $T_{I}$ action on  
$\C^{|I|}$. Let us consider 

\[ V_{I} \times \C^{|I|},\]
on which a torus action $\T_{I}$ acts diagonally.  Let  $\chi_{I}$ be a map defined by:

\[ \chi_{I} : V_{I} \times \C^{|I|} \to \R^{|I|} : \chi_{I}(m, t, z) = (\alpha_{i} - t\varphi_{i}(m) - |z_{i}|^{2})_{i \in I}.  \]
 It is clear that 0 is a regular value for $\chi_{I}$ and $\T_{I}$ is locally free on  $\chi_{I}^{-1}(0)$. Therefore, 
 
\begin{equation}
\label{a11}
W_{I} = \chi_{I}^{-1}(0) / \T_{I}
\end{equation}
defines an open orbifold with boundary isomorphic to $U_{I}$.

For the second part of the lemma, let $V_{I}^{J}$ be an orbifold defined by 

\[ V_{I}^{J} =  \{ (m, t) \in U_{J} \times [1, \infty) \big| t \cdot \varphi_{i}(m) < \alpha_{i}, \mathrm{for \ all} \ i \in \check{I} \}. \]
As the construction of $W_{J}$, we can define a map $\tilde{\chi}_{J}$ on 
$V_{I}^{J} \times \C^{|J|}$ and 

\[ W_{I}^{J} = \tilde{\chi}_{J}^{-1}(0) / \T_{J}, \]
gives an orbifold with boundary isomorphic to $U_{J}$. 

For $I \subsetneq J$, we have that $U_{I} \subseteq U_{J}$ and $\check{J} \subseteq \check{I}$. Hence, $V_{I}^{J}$ is an open subset of both $V_{J}$ and $V_{I}$.  And  $W_{I}^{J}$ can be identified as open suborbifold of both $W_{I}$ and $W_{J}$.  We define the diffeomorphism  $\Phi_{I}^{J}$ to be the map  between $W_{I}$ and $W_{J}$ factoring through $W_{I}^{J}$. The verification of (\ref{a13}) is straightforward. 
\end{proof}

\end{lemma}

Thanks to Lemma \ref{a14}, we can obtain an orbifold $W$ by gluing all the $\{ W_{I} \}_{I \in \mathcal{A}}$ using $\phi_{I}^{J}$. From the construction, one can check that $W$ is a compact orbifold with boundary isomorphic to $\Sigma$. Therefore, we have the following theorem. 

\begin{theorem}

\label{a16}
Let $\Sigma$ be a compact orbifold. If it has a good cover, then we can construct an orbifold  $W$ with boundary isomorphic to $\Sigma$. 

\end{theorem}

Now, let $(M, E, \phi)$ be a $K$-cycle with compact vanishing set and $\Sigma = \partial M$.   We will show that $\Sigma$ has a good cover. 

\begin{definition}
Let $H$ be an isotropy group of $\T$-action on $\Sigma$ and

\[ \Sigma_{H} = \{x \in \Sigma \big| \T_{x} = H \}, \]
where $\T_{x}$ is the isotropy group of $x$. For any connected component $F$ of  $\Sigma_{H}$, we denote by $U_{F}$ a $\T$-invariant neighborhood of $F$ in $\Sigma$. Here, we can choose $U_{F}$ small enough  such that  

\begin{itemize}
\item 
For all $x \in U_{F}$, we have that $\T_{x} \subseteq H$. 
\item
 For any other component $F^{'}$ of $\Sigma^{H}$, we have that $U_{F} \cap U_{F^{'}} = \emptyset$. 
\end{itemize}
In this case, we say that $U_{F}$ has $\it{level}$ equal to $\mathrm{dim}(H)$. As $H$ ranges over all the isotropy group, we obtain  an open cover of $\Sigma$, denoted by $\{U_{i} \}$. 

\end{definition}

\begin{lemma}
\label{a3}
Every point $x \in \Sigma$ can  be covered by at most $\mathrm{dim}(\T) - \mathrm{dim}(\T_{x})$ open sets in $\{ U_{i} \}$. 
\begin{proof}
The lemma follows from the fact that every point $x \in \Sigma$ can only be covered by open set $U_{i}$ with level no lower than $\mathrm{dim}(\T_{x})$.
\end{proof}

\end{lemma}

\begin{proposition}
\label{a15}
There exists a good cover on $\Sigma$.  
\begin{proof}
It is enough to associate every $U_{i}$ with a compatible circle action. We will complete the proof by induction on the level of $\{U_{i}\}$. 

If $U_{i}$ has the highest level, then we can always find a circle action $S_{i}$ which acts locally freely on $U_{i}$, using the fact that $\phi$ induces a nowhere vanishing vector field on $\Sigma$. 

Suppose that we have already associated the open sets whose level is greater than $K$ with compatible circle actions.  
 
Let $U_{k}$ be an open set  with level $K$. Suppose that $I$ is a subset of $\{1, \dots, n\}$ such that  

 \[ U_{k} \cap ( \bigcap_{i \in I} U_{i}) \neq \emptyset, \]
 and every $\{ U_{i} \}_{i \in I}$ has level greater than $K$. From Lemma \ref{a3},  we know that $|I| \leq \mathrm{dim}(\T) -1- K$. Meanwhile, for any point $x \in U_{k}$, the isotropy group $T_{x}$ has dimension no greater than $K$. Hence, the compatible circle action $S_{k}$ always exists. This completes the proof. 
\end{proof}
\end{proposition} 

In order to  construct a cap for the $K$-cycle $(M, E, \phi)$, we will show how to define an orbifold vector bundle $L$ and  an equivariant map $\psi$ on $W$. 

\begin{proposition}
\label{a4}
Fixed any $\T$-weight $\gamma$, we can construct a $K$-cycle $(W, L, \psi)$ such that 

\[ \partial (W, L, \psi) \cong (\Sigma, E, \phi)\]
and $L|_{x}$ has $\T$-weight equal to $\gamma$ for any $x \in W^{\T}$. 

\begin{proof}
According to the construction in Theorem \ref{a16}, it is enough to define the vector bundle and equivariant map on every piece $W_{I}$. 

Let $\gamma_{I}$ be the  restriction of $\gamma$ to $S_{I}$, which is subgroup of $\T$. Let $F_{\gamma_{I}}$  be the trivial line over $\C^{|I|}$ on which $T_{I}$ acts with weight $\gamma_{I}$. Recall that 
$\chi_{I}^{-1}(0)$ is an open subset of $V_{I} \times \C^{I}$ and the diagonal action $\T_{I}$ is locally free on $\chi_{I}^{-1}(0)$.  Thus,

\begin{equation}
\label{t1}
L_{I} = ((E \boxtimes F_{\gamma_{I}})|_{\chi_{I}^{-1}(0)}) / \T_{I}
\end{equation}
defines an orbifold vector bundle on $W_{I} = \chi_{I}^{-1}(0)/ \T_{I}$.  

For the equivariant map, let $\hat{\phi}$ be the pullback of $\phi$ to $V_{I} \times \C^{|I|}$. After restricting to $\chi_{I}^{-1}(0)$, $\hat{\phi}$ descends to a map on $W_{I}$. Using the diffeomorphism in Lemma \ref{a14}, we can get an orbifold vector $L$ and an equivariant map $\psi$ by gluing. And the triple $(W, L, \psi)$ constitute a $K$-cycle satisfying all the desired properties. 
\end{proof}

\end{proposition}

Next, we are going to show that $(W, L, \psi)$ is bordant to a discrete $K$-cycle. 

\begin{definition}
Let $\D^{n} = D^{2} \times \dots \times D^{2}$ be the product of $n$ copies of open disks. A $n$-dimensional torus $\T^{n} = S^{1} \times \dots \times S^{1}$ acts on $\D^{n}$ in such a way that the $i$-th factor of $S^{1}$ acts on the $i$-th Disk by rotation.  
\end{definition}

For every $I \in \mathcal{A}$,  let us define $K_{I} = \rho_{I}^{-1}(\alpha_{I}) \cap V_{I}$, that is 

\[ K_{I} = \{ (m,t) \in U_{I} \times [1, \infty) \big| t\cdot \varphi_{i}(m) = \alpha_{i}, i \in I \,  \mathrm{and} \  t\cdot \varphi_{j}(m) < \alpha_{j}, j \in \check{I} \}. \]
When $\check{I} = \emptyset$, $K_{I}$ is a compact orbifold without boundary. Otherwise, it is possible that $\partial K_{I} \neq \emptyset$. However, we have the following. 

\begin{lemma}
\label{a24}
For every $I \in \mathcal{A}$, we can obtain a closed orbifold by gluing 

\[  (\bigsqcup_{J \in \mathcal{A}, I \subsetneq J} K_{J} \times_{\T^{I}_{J}} \D^{|J|-|I|}) \sqcup K_{I}, \]
where $\T^{I}_{J}$ is a $(|J|-|I|)$-dimensional torus and it acts on $K_{J}$ through $S_{J}^{I} = \prod_{i \in J \cap \check{I}}S_{i}$. 
\begin{proof}
It is straightforward based on the fact that $\alpha_{I}$ are regular values for 

\[ \rho_{I} : U_{I} \times [1, \infty) \to \R^{|I|} .\]
\end{proof}
 
\end{lemma}

\begin{lemma}
\label{a26}
The $K$-cycle $(W, L, \psi)$ is bordant to a $K$-cycle in the following form:

\[ (W, L, \psi) \sim \sum (M_{i} \times_{\T_{i}}\D^{|\T_{i}|}, E_{i}, \phi_{i}) + \sum (N_{j}, F_{j}, \psi_{j}),\] 
where $\{ M_{i}, N_{j} \}$ are compact oribifolds without boundary. 

\begin{proof}
If we identify $K_{I}/\T_{I}$ with a subset in $W_{I}$,  then the set 
\begin{equation}
\label{a23}
 Z_{I} = K_{I} \times_{\T_{I}} \D^{|I|}
\end{equation}
can be identified as  a neighborhood of $K_{I}/ \T_{I}$ in $W_{I}$.  We define $Z$ to be the orbifold obtained by gluing $Z_{I}$ as in Lemma \ref{a14}.  Since $V^{\phi}$ is nowhere vanishing on $\Sigma$, the vanishing set $V^{\psi}$ in $W$ must be contained in $Z$. Therefore, 

\[ (W, L, \psi) \sim (Z, L|_{Z}, \psi|_{Z}). \]
Without losing generality, we can furthermore assume that 

\begin{equation}
\label{a25}
(Z, L|_{Z}, \psi|_{Z}) \sim \sum_{I}(Z_{I}, L|_{Z_{I}}, \psi|_{Z_{I}}). 
\end{equation}
Therefore, the Lemma follows from (\ref{a23}), (\ref{a25}), and Lemma \ref{a24}.

\end{proof}

\end{lemma}

By repeatedly using Lemma \ref{lemma-2}, the following  theorem follows from Lemma \ref{a26}. 

\begin{theorem}
The $K$-cycle $(W, L, \psi)$ is bordant to a discrete $K$-cycle. 

\end{theorem}


\subsection{Nonabelian Case}
Now, we assume that $G$ is a compact connected Lie group, $\T$ is a maximal torus, and  $\mathfrak{t}_{+}$ is a fixed positive Weyl chamber.  For any $x \in \mathfrak{t}_{+}$, we denote by $G_{x}$ the isotropy group of adjoint action.  

\begin{lemma}
Let $\Delta$ be a face of $\mathfrak{t}_{+}$. For any point in $\mathrm{int}(\Delta)$, they have the same isotropy group, denoted by $G_{\Delta}$. Moreover, we have that $G_{\Delta}  \subset G_{\Delta^{'}}$ if and only if $\Delta$ is a sub-face of $\Delta^{'}$.
\end{lemma}

Let us recall symplectic cross-section theorem.

\begin{theorem}[Cross-section]
Let $(M, \omega)$ be a compact connected symplectic orbifold  with a moment map $\mu : M \to \mathfrak{g}^{*}$ arising from an action of a compact Lie group G. For any face $\Delta$ in $\mathfrak{t}_{+}$, let $V_{\Delta}$ be a small neighborhood of $\mathrm{int}(\Delta)$ in $\mathfrak{t}_{+}$ such that $G_{x} \subseteq G_{\Delta}$ for any $x \in V_{\Delta}$. If we denote $U_{\Delta} = G_{\Delta} \cdot V_{\Delta}$, then the cross section $R = \mu^{-1}(U_{\Delta})$ is a $G_{\Delta}$-invariant symplectic sub-orbifold and 

\[ U = G\cdot R = G \times_{G_{\Delta}} R\] 
is an open subset of $M$. Moreover, if $A_{\Delta}$ is the abelian part of $G_{\Delta}$, then the $A_{\Delta}$-action on $R$ extends in a unique way to an action on $U$ which commutes with the G-action.
\begin{proof}
See \cite{Guillemin90} \cite{Lerman98}. 
\end{proof}
\end{theorem}

In this paper, we are considering stably complex orbifolds instead of symplectic orbifolds. Hence, the symplectic cross-section theorem does not apply. However, the idea of building the cap in nonabelian case comes from the symplectic cross-section theorem and symplectic surgery by Meinrenken \cite{Meinrenken98}.  

Let $(M, E, \phi)$ be a $K$-cycle with compact vanishing set and $\Sigma = \partial M$. 
\begin{definition}
For each $m \in \Sigma$, let $\mathfrak{g}_{m} \subset \mathfrak{g}$ be the corresponding isotropy Lie algebra. It is clear that $\mathfrak{g}_{g \cdot m} = \mathrm{Ad}(g)(\mathfrak{g}_{m})$. We call the set of subalgebras

\[ (\mathfrak{g}_{m}) = \{ \mathrm{Ad}(g)(\mathfrak{g}_{m}) \big| g \in G\} \]
the \it{orbit  type} of $m$. There are only finite many orbit type on $\Sigma$. Moreover, there is a unique orbit type $(\mathfrak{g}_{0})$ such that the set 

\[ \Sigma_{(\mathfrak{g}_{0})} = \{ m \in \Sigma \big|(\mathfrak{g}_{m}) = (\mathfrak{g}_{0}) \} \] 
is a dense, open subset in $\Sigma$ \cite{Guillemin90}\cite{Lerman98}. 
\end{definition}

  Let $(\mathfrak{h})$ be an arbitrary orbit type of $\Sigma$. There exists a face $\Delta$ with maximum dimension such that $(\mathfrak{h})$ is subconjugated to  $(\mathfrak{g}_{\Delta})$.   Suppose that $F$ is a connected component of $\Sigma_{(\mathfrak{h})}$ and $U$ is a small $G$-invariant neighborhood of $F$ in $\Sigma$ such that for any $m \in U$
  
\[ (\mathfrak{g}_{m}) \subseteq (\mathfrak{g}_{\Delta}).\]  
In this case, as in cross-section theorem, we can find a $G_{\Delta}$-invariant subset $R$ of $U$ such that 
  
  \[   U =  G \cdot R \cong G \times_{G_{\Delta}} R, \]  
where the isomorphism map is given by 

\[ G \times_{G_{\Delta}} R \to G \cdot R, \ [a, u] \to a \cdot u. \]
Moreover, the $A_{\Delta}$-action on $R$ extends  to an action on $U$, which commutes with the $G$-action \cite{Woodward96} \cite{Meinrenken98}.  

\begin{lemma}
\label{a27}
There exists an open coving $\{ U_{i} \}$ of  $\Sigma$ with circle actions $\{S_{i} \}$ such that  
\begin{itemize}
\item
    Every $U_{i}$ is $G$-invariant. 
\item    
   The circle action $S_{i}$ acts locally freely on $U_{i}$ and commutes with $G$-action.
\item
  For all $I \in \mathcal{A}$, $\{ S_{i} \}_{i \in I}$ generate an $|I|$-dimensional torus and it acts locally freely on $U_{I}$. 
\end{itemize}
\begin{proof}
The proof is similar to the torus case except the circle actions are induced from the locally defined action $A_{\Delta}$. 
\end{proof}

\end{lemma}

By similar argument in the torus case,  we can prove the following theorem.

\begin{theorem}
There exists a $K$-cycle $(W, L, \psi)$ such that 

\[  \partial(W, L, \psi) \cong (\Sigma, E, \phi). \]
Moreover, $(W, L, \psi)$ is bordant to a discrete $K$-cycle. 
\end{theorem}


\section{Quantization Map}
In this section, we will first give the definition of quantization map for  all $K$-cycles:

$$  Q_{\mathrm{TOP}}  :   \{(M, E, \phi)\}  \To \hat{R}(G).$$
Then, we will  show that $Q_{\mathrm{TOP}}$ induces  an isomorphism from $\hat{K}(G)$ to $\hat{R}(G)$. 

To begin with, by Theorem \ref{localization theorem}, we know that  

\begin{equation}
\label{t11}
(M, E, \phi)  \sim  \bigsqcup_{k} (U_{k}, E|_{U_{k}}, \phi|_{U_{k}}), 
\end{equation}
where the right-hand side consists of $K$-cycles with compact vanishing set.  For each $(U_{k}, E|_{U_{k}}, \phi|_{U_{k}})$, we can build a cap $(W_{k}, L_{k}, \psi_{k})$. Recall that the constructions of $L_{k}$ are not unique. In fact, we can build $L_{k}$ in a way such that $\{ (W_{k}, L_{k}, \psi_{k}) \}$ forms a global $K$-cycle by putting all them together.

\begin{lemma}
\label{t14}
We can build $(W_{k}, L_{k}, \psi_{k})$ such that

\[ \bigsqcup_{k} (W_{k}, L_{k}, \psi_{k}) \]
is a $K$-cycle. 
\begin{proof}
Let us denote $\Sigma_{k} = \partial U_{k}$. Every $W_{k}$ is constructed as in section 4 and $\psi_{k}$ is induced from $\phi|_{\Sigma_{k}}$. Since $\phi$ is proper over $M^{\phi}$, we can choose $U_{k}$ small enough so that  the maps $\{ \psi_{k} \}$ are proper over $W = \bigsqcup_{k} W_{k}$. The problem left is to construct line bundles $\{ L_{k} \}$ which are compatible with $\{ \psi_{k} \}$ globally. That is, there exists a constant $C$ such that for all 
$k$,

\begin{equation}
\label{t2}
 \|  \frac{\sqrt{-1}}{2\pi}L_{\xi}  + \langle \psi_{k}(m) , \xi \rangle \| \leq C \cdot \| \xi \|,   \ m \in W_{k}^{\psi_{k}} \ \text{and} \ \xi \in \mathfrak{g}_{m}. 
 \end{equation}

Due to Remark \ref{s30}, there exists a constant $R$ and a series of dominant weights $\{ \gamma_{k} \} \in \mathfrak{t}_{+}$ such that 

\begin{equation}
\label{t3}
\phi_{k}(\Sigma_{k}) \cap \mathfrak{t}_{+} \subseteq  B(\gamma_{k}, R) = \{ x \in \mathfrak{t}_{+} \big| \| x - \gamma_{k} \| \leq R \} 
\end{equation}
and $\lim_{k \to \infty} \| \gamma_{k} \| = \infty$. As in Proposition \ref{a4}, we can construct $L_{k}$  using the fixed weights $\gamma_{k}$. These $\{ L_{k} \}$ satisfy condition (\ref{t2}).

\end{proof}
\end{lemma}

From section 4, we know that every $(W_{k}, L_{k}, \psi_{k})$ is bordant to a discrete $K$-cycle:

\[ (W_{k}, L_{k}, \psi_{k}) \sim \sum_{i} (N_{k}^{i}, F_{k}^{i}, \rho_{k}^{i}). \]
It is natural to ask that if we put all them together, do they form a $K$-cycle?

\begin{lemma}
\label{a33}
The infinite sum 

\[ \sum_{k, i} (N_{k}^{i}, F_{k}^{i}, \rho_{k}^{i}) \]
constitute a $K$-cycle. 
\begin{proof}
Let $N = \bigsqcup_{k, i} N_{k}^{i}$,  $F$ be the orbifold vector bundle such that $F|_{N_{k}^{i}} = F_{k}^{i}$, and $\rho$ be the map on $N$ such that $\rho|_{N_{k}^{i}} = \rho_{k}^{i}$. To check that $(N, F, \rho)$ is a $K$-cycle, it is enough to show that $\rho$ is proper. To be more precise, for every $R>0$, we need to show that there exists constant $T_{R}$ such that  for all $i, k \geq T_{R}$,

\[ \rho_{k}^{i} (N_{k}^{i}) \cap B_{R} = \emptyset,  \]
 where 

\[ B_{R} = \{  x \in \mathfrak{t}_{+} \big| \|x \| \leq R \}. \]

First,  given any fixed $k$, we know that

\[ \rho_{k}^{i} (N_{k}^{i}) \cap B_{R} = \emptyset,  \]
when $i$ is large enough. On the other hand, there exists a constant $R^{'} \gg R$ such that if  $\psi_{k}(W_{k}) \cap B_{R^{'}} = \emptyset$,  then 

\[ \rho_{k}^{i} (N_{k}^{i}) \cap B_{R} = \emptyset  \ \mathrm{for \ all} \ i. \]
By  (\ref{t3}), we have that $\psi_{k}(W_{k}) \cap B_{R^{'}} = \emptyset$ as $k \to \infty$. This completes the proof.

\end{proof}

\end{lemma}

\begin{theorem}
\label{t8}
Every $K$-cycle $(M, E, \phi)$ is bordant to a discrete $K$-cycle:

\begin{equation}
\label{k4}
 (M, E, \phi) \sim \sum_{k=1}^{\infty} (M_{k}, E_{k}, \phi_{k})
 \end{equation}
In addition, for any irreducible representation $\gamma \in \hat{G}$, there exists a constant $T_{\gamma}$ such that 

\begin{equation}
\label{o10}
[\mathrm{Index}(M_{k}, E_{k})]^{\gamma} = 0, \ \mathrm{if} \  k \geq T_{\gamma}. 
\end{equation}

\begin{proof}
The first part follows from Lemma \ref{a33}. For the section part, it is enough to prove the abelian case. The general case follows from the induction argument in \cite{Paradan01}. 

Suppose that $G = \T$ is abelian. By the discussion before, there exists a constant $R$ and a series of weights $\{ \gamma_{k} \} \in \mathfrak{t}$ such that 

\begin{equation}
\label{t19}
\phi_{k}(M_{k}) \subseteq  B(\gamma_{k}, R) = \{ x \in \mathfrak{t} \big| \| x - \gamma_{k} \| \leq R \} 
\end{equation}
and $\lim_{k \to \infty} \| \gamma_{k} \| = \infty$. Therefore,  we know that when $k$ is large enough, there exists a cyclic unit vector $\xi \in \mathrm{Lie}(\T)$ such that 

\[ \langle \phi_{k}(x), \xi \rangle \gg 0, \ \mathrm{for \ all} \ x \in M_{k}^{\T}. \]
By (\ref{s27}), we conclude that as $k \to \infty$

\[ \langle \omega^{k}_{x}, \xi \rangle \gg 0, \ \mathrm{for \ all} \ x \in M_{k}^{\T}, \]
where $\omega^{k}_{x}$ is the fiber weight of $E_{k}|_{x}$. Therefore, we can finish the proof by the Atiyah-Bott fixed point theorem (for details, one can see Theorem 5.1 in \cite{Meinrenken99}).

\end{proof}
\end{theorem}

Now, we are able to  define the quantization map. 

\begin{definition}
\label{t16}
Let $\bigsqcup_{k=1}^{\infty} (M_{k}, E_{k}, \phi_{k})$ be a discrete $K$-cycle satisfying (\ref{o10}),  we define the quantization map $Q_{\mathrm{TOP}}$ to be

\[ Q_{\mathrm{TOP}}\big( \bigsqcup_{k=1}^{\infty} (M_{k}, E_{k}, \phi_{k})\big) = \sum_{k} \mathrm{Index}(M_{k}, E_{k}) \in \hat{R}(G).\]

\end{definition}

Suppose that $(M, E, \phi)$ is an arbitrary $K$-cycle. As in (\ref{k4}), let us assume that 

\begin{equation}
\label{t12}
(M, E, \phi) \sim \sum_{k}(M_{k}, E_{k}, \phi_{k}).
\end{equation}
If we can prove that $Q_{\mathrm{TOP}}$ is invariant under bordism, then we can get a well-defined quantization map for any $K$-cycle:

\[ Q_{\mathrm{TOP}}(M, E, \phi) = Q_{\mathrm{TOP}}\big( \bigsqcup_{k=1}^{\infty} (M_{k}, E_{k}, \phi_{k})\big) \]



\begin{lemma}
\label{t22}
Suppose that $(W, L, \psi)$ is a $K$-chain. For any constant $R > 0$, there exists a hypersurface $\Sigma_{R}$ in W such that 

\begin{itemize}
\item 
  The vector field $V^{\psi}$ is nowhere vanishing over $\Sigma_{R}$;
  
\item  
$\Sigma_{R}$ subdivides W into two parts: a bounded part $W_{-}$ and an unbounded part $W_{+}$,  with the property  that 

\[ \psi(W_{+}) \cap B_{R}  = \emptyset . \]   
\end{itemize}
\begin{proof}
Let $\mathcal{H}$ be a function defined by $\mathcal{H} = \| \psi \|^{2} : W \to \R$. By Theorem \ref{separate theorem},  the vanishing set $W^{\psi}$ can be covered by $\{ U_{\alpha} \}$ such that 
$F_{\alpha}= W^{\psi} \cap U_{\alpha}$ is compact. Further, by Remark \ref{s30},  there exists a constant $K$ such that for all $\alpha$,

\begin{equation}
\label{t21}
  \|  \psi(x) - \psi(y) \| \leq K, \mathrm{for \ all} \ x, y \in U_{\alpha} \cap \psi^{-1}(\mathfrak{t}_{+}). 
\end{equation}
By (\ref{t21}) and the fact that $\psi$ is proper over the vanishing set, there are only finitely many  $\alpha$ such that 

\begin{equation}
\label{q1}
 F_{\alpha} \cap \mathcal{H}^{-1}(R + K) \neq \emptyset.
 \end{equation}
Let us denote them by $\{ \alpha_{1}, \dots, \alpha_{n} \}$. Define a G-invariant non-negative function $\rho : W \to \R$ such that 

\begin{itemize}
\item
  $\rho(x) = 0$ for all $x \notin \bigsqcup_{i=1}^{n} U_{\alpha_{i}}$ 
  
\item 
  $ \rho(x) \geq 2(R+K)$ for  $x \in U_{\alpha_{i}} \cap W^{\psi}$, $i = 1, \dots, n$. 
\end{itemize}
Select a regular value $c$ for $\mathcal{H} + \rho$, which is very close to $R+K$. The hypersurface $\Sigma_{R} = (\mathcal{H} + \rho)^{-1}(c)$ satisfies all the conditions. 

\end{proof}

\end{lemma}

\begin{proposition}
Let  $(M, E, \phi)$ and $(M^{'}, E^{'}, \phi^{'})$ be two $K$-cycles. If they are bordant, then

$$ Q_{\mathrm{TOP}}(M, E, \phi) = Q_{\mathrm{TOP}}(M^{'}, E^{'}, \phi^{'}) \in \hat{R}(G).$$

\begin{proof} 
Fixing an irreducible representation $\gamma \in \widehat{G}$, it is enough to show that 

\[ [Q_{\mathrm{TOP}}(M, E, \phi)]^{\gamma} = [Q_{\mathrm{TOP}}(M^{'}, E^{'}, \phi^{'})]^{\gamma}.\]  
By Definition \ref{t12}, let us assume that 

\[ Q_{\mathrm{TOP}}(M, E, \phi) = \sum_{k=1}^{\infty} \mathrm{Index}(M_{k}, E_{k}) \]
and 

\[ Q_{\mathrm{TOP}}(M^{'}, E^{'}, \phi^{'}) = \sum_{k=1}^{\infty} \mathrm{Index}(M_{k}^{'}, E_{k}^{'}). \]
Since $(M, E, \phi)$ is bordant to $(M^{'}, E^{'}, \phi^{'})$, we also have 

\begin{equation}
\label{o5}
\bigsqcup_{k}(M_{k}, E_{k}, \phi_{k}) \sim \bigsqcup_{k}(M_{k}^{'}, E_{k}^{'}, \phi_{k}^{'}). 
\end{equation}
Suppose $(W, L, \psi)$ gives the bordism in (\ref{o5}).  

As in Lemma \ref{t22}, we can construct a hypersurface $\Sigma_{R}$ which cuts $(W, L, \psi)$ into two pieces: $(W_{-}, L|_{W_{-}}, \psi|_{W_{-}})$ and $(W_{+}, L|_{W_{+}}, \psi|_{W_{+}})$.  We observe that $W_{-}$ is an orbifold with boundary consisting of 

\begin{itemize}
\item 
 The part which doesn't intersect with $\Sigma_{R}$:  $\bigsqcup_{k \in \mathscr{A}} M_{k} $ and $\bigsqcup_{j \in \mathscr{A}^{'}}M_{j}^{'}$.  
\item
 The part which intersects with $\Sigma_{R}$ : $\bigsqcup_{k \in \mathscr{B}} U_{k} $ and $\bigsqcup_{j \in \mathscr{B}^{'}}U_{j}^{'}$, where $U_{k}, U_{j}^{'}$ are subsets of   $M_{k}$ and  $M_{j}^{'}$.  
\end{itemize}
For any fixed $\gamma \in \widehat{G}$, we can choose $R$ large enough such that for all $k \notin \mathscr{A}$ and $j \notin \mathscr{A}^{'}$

\[ [Q(M_{k}, E_{k})]^{\gamma} = [Q(M_{j}^{'}, E_{j}^{'})]^{\gamma} = 0. \]

Notice that $(W_{-}, L|_{W_{-}}, \psi|_{W_{-}})$ is a $K$-chain with compact vanishing set.  As the construction in Section 4,  we can build a cap and obtain a compact $K$-chain (which gives a compact bordism) by gluing on the cap. During the gluing process, we notice that 

\[ \{ (M_{k}, E_{k}, \phi_{k}) \}_{k \in \mathscr{A}} \ \mathrm{and} \{ (M_{j}^{'}, E_{j}^{'}, \phi_{j}^{'}) \}_{j \in \mathscr{A}^{'}} \]
remain the same. When the constant $R$ is large enough, we can also assume that the multiplicity of $\gamma$ does not change during the gluing. Therefore, due to the fact that index map is invariant under compact bordism, we can conclude that 

\[ [Q_{\mathrm{TOP}}(M, E, \phi)]^{\gamma}  = [Q_{\mathrm{TOP}}(M^{'}, E^{'}, \phi^{'})]^{\gamma} \in \Z. \]
\end{proof}

\end{proposition}

\begin{proposition}
Given any $K$-cycle $(M, E, \phi)$, we have 

$$ Q_{\mathrm{TOP}}(M, E, \phi) = Q_{\mathrm{TOP}}(\hat{M}, \hat{E}, \hat{\phi}), $$
where $(\hat{M}, \hat{E}, \hat{\phi})$ is the bundle modification in Definition \ref{d1}. 

\begin{proof}
Assume that $M$ is modified by the principle bundle $P$ with fiber $N$. First, by Definition \ref{t16}, we have

\[ (M, E, \phi)  \sim \sum(M_{k}, E_{k}, \phi_{k}), \]
Suppose that $(W, L, \psi)$ gives the bordism above. As  the construction of  vector bundle  over $W$, we can extend the principle bundle $P$ on $M$ to $W$, denoted by $P_{W}$. Now, we can define the bundle modification (with respect to $P_{W}$) of $W$ to be

\[ \hat{W} =  P_{W} \times_{H} N .\]  
Let $\hat{L}$ and $\hat{\psi}$ be the pullback of $L$ and $\psi$ to $\hat{W}$. Then $(\hat{W}, \hat{L}, \hat{\psi})$ is a $K$-chain, giving a bordism between 

\[ (\hat{M}, \hat{E}, \hat{\phi}) \ \text{and} \  \sum_{k} (\hat{M}_{k}, \hat{E_{k}}, \hat{\phi_{k}}).  \]
By the definition of quantization, we have 

\[ Q_{\mathrm{TOP}}(\hat{M}, \hat{E}, \hat{\phi} ) =  \sum_{k} \mathrm{Index}(\hat{M}_{k}, \hat{E_{k}}).  \]
On the other hand, because $(\hat{M}_{k}, \hat{E}_{k})$ is bundle modification of $(M_{k},E_{k})$,  we have that

\[ \mathrm{Index}(M_{k}, E_{k})  = \mathrm{Index}(\hat{M}_{k}, \hat{E_{k}}).\] 
This completes the proof. 
\end{proof}

\end{proposition}

By the two propositions above, we can conclude that the quantization map $Q_{\mathrm{TOP}}$ is a well-defined map from $\hat{K}(G)$ to $\hat{R}(G)$.

\begin{theorem}
The quantization map $Q_{\mathrm{TOP}}$ gives an isomorphism:

\[ Q_{\mathrm{TOP}}  : \hat{K}(G) \cong \hat{R}(G). \]
\begin{proof} Let $P$ be the map defined in Theorem \ref{o9}. It is clear that $Q_{\mathrm{TOP}}  \circ P : \hat{R}(G) \To \hat{R}(G)$ is the identity map. Therefore, $Q_{\mathrm{TOP}}$ is surjective.  

Since every $K$-cycle is equivariant to a discrete $K$-cycle, it is enough to prove the injectivity of compact $K$-cycles. This immediately follows from the fact that  geometric $K$-homology is isomorphic to analytic $K$-homology \cite{Baum10} \cite{Baum07}. 

\end{proof}
\end{theorem}

\section{Quantization Commutes with Reduction}
In the previous section, we give the definition of quantization for general $K$-cycles. In this framework, we will provide a new approach to the Quantization Commutes with Reduction theorem in non-compact setting.  Let us begin with the multiplicative property. 

\begin{theorem}
\label{t17}
Let $(M,E, \phi_{1})$ and $(N, F, \phi_{2})$ be two K-cycles in $\hat{K}(G)$, where $N$ is compact. We have the following: 

$$ Q_{\mathrm{TOP}}(M, E, \phi_{1}) \times Q_{\mathrm{TOP}}(N, F, \phi_{2}) = Q_{\mathrm{TOP}}(M \times N, E \boxtimes F, \phi_{1} \circ \pi_{1} + \phi_{2} \circ \pi_{2}) $$
where $\pi_{1}$ is the projection from $M \times N$ to M and where $\pi_{2}$ is the one to N. 

 \begin{proof}
Consider the $K$-chain $(M \times N \times [0, 1] , \hat{E} \otimes \hat{F}, \phi)$, where $\hat{E}$ and $\hat{F}$ are the pullback of line bundles $E$ and $F$, and the map $\phi$ is defined by 

$$ \phi(m, n, t) = \phi_{1}(m) + t\phi_{2}(n). $$
This gives a bordism between

\begin{equation}
\label{k2}
(M \times N, E \boxtimes F, \phi_{1} \circ \pi_{1}) \ \text{and} \ (M \times N, E \boxtimes F, \phi_{1} \circ \pi_{1} + \phi_{2} \circ \pi_{2}).
\end{equation}
By Theorem \ref{t8}, let us assume that there is a bordism 

\begin{equation}
\label{o6}
(M, E, \phi_{1}) \sim  \bigsqcup_{k} (M_{k}, E_{k}, \rho_{k}). 
\end{equation}
Since $N$ is compact, (\ref{o6}) induces another bordism:

\[ (M \times N, E \boxtimes F, \phi_{1} \circ \pi_{1}) \sim  \bigsqcup_{k} (M_{k} \times N, E_{k} \boxtimes F, \hat{\rho_{k}}), \]
where $\hat{\rho_{k}}$ is the composition of $\rho_{k}$ with the projection from $M_{k} \times N$ to $M_{k}$. By the bordism invariance of quantization, we have 

\[ Q_{\mathrm{TOP}}(M \times N, E \boxtimes F, \phi_{1} \circ \pi_{1}) = \sum_{k} \mathrm{Index}(M_{k} \times N, E_{k} \boxtimes F).  \]
Meanwhile, we have that 

\[ \sum_{k} \mathrm{Index}(M_{k} \times N, E_{k} \boxtimes F) = ( \sum_{k} \mathrm{Index}(M_{k}, E_{k})) \times \mathrm{Index}(N, F). \]
Hence, we get

\begin{equation}
\label{k1}
Q_{\mathrm{TOP}}(M \times N, E \boxtimes F, \phi_{1} \circ \pi_{1}) = Q_{\mathrm{TOP}}(M, E, \phi_{1}) \times Q_{\mathrm{TOP}}(N, F, \phi_{2}).
\end{equation}
The theorem follows from (\ref{k2}) and (\ref{k1}).
 
\end{proof}

\end{theorem}

\begin{theorem}
Let $(M, E, \mu)$ be pre-quantum data. If the moment map $\mu$ is proper and 0 is a regular value, then 

\[ [Q_{\mathrm{TOP}}(M, E, \mu)]^{G} = Q_{\mathrm{TOP}}(M_{0}, E_{0}),  \]
where $M_{0} = \mu^{-1}(0)/G$ and $E_{0} = (E|_{\mu^{-1}(0)})/G$. 
\end{theorem}

\noindent  $\mathbf{Sketch \ of \ the  \ proof}$:
First, recall the decomposition of the vanishing set $M^{\mu}$  \cite{Paradan01}:
 
 \[ M^{\mu} = \sum_{\gamma \in \Gamma} G . (M^{\gamma} \cap \mu^{-1}(\gamma)), \] 
where $\Gamma$ is a discrete set in $\mathfrak{t}_{+}$. For each $\gamma \in \Gamma$, let  $U_{\gamma}$ be a small open G-invariant  neighborhood of $G . (M^{\gamma} \cap \mu^{-1}(\gamma))$ in $M$. 
By Theorem \ref{localization theorem}, we have 

\[ Q_{\mathrm{TOP}}(M, E, \mu) = \sum_{\gamma \in \Gamma} Q_{\mathrm{TOP}}(U_{\gamma}, E|_{U_{\gamma}}, \mu|_{U_{\gamma}}). \]
In accord with \cite{Paradan01}, one can show  that for all $\gamma \neq 0$, 

\[ [Q_{\mathrm{TOP}}(U_{\gamma}, E|_{U_{\gamma}}, \mu|_{U_{\gamma}})]^{G} = 0. \]
Thus,

\[ [Q_{\mathrm{TOP}}(M, E, \mu)]^{G} = [Q_{\mathrm{TOP}}(U_{0}, E|_{U_{0}}, \mu|_{U_{0}})]^{G}. \]
Since $0$ is a regular value, we have that $U_{0} \cong \mu^{-1}(0) \times \mathfrak{g}^{*}$. As in Section 4, we can build a cap  $(W, L, \psi)$. By gluing on the cap, we can get a compact $K$-cycle $(N, F)$. Remember that the construction of the cap (in particular, the line bundle $L$) is not unique.  In fact, we can build $(W, L, \psi)$ in a way such that 

\[ \mathrm{Index}(N, F) = [Q_{\mathrm{TOP}}(M, E, \mu)]^{G}. \]
In this case,  we can show that $(N, F)$ is in fact a bundle modification of $(M_{0}, E_{0})$ (for example, when $G = S^{1}$, $N$ is a sphere bundle over $M_{0}$). Therefore, 

\[Q_{\mathrm{TOP}}(M_{0}, E_{0}) =[Q_{\mathrm{TOP}}(M, E, \mu)]^{G} .\]

As a result, Theorem \ref{o7} follows from the two theorems above.



\bibliographystyle{alpha}
\bibliography{mybib}

\end{document}